\documentclass[a4paper,11pt,twoside,reqno]{amsart}

\usepackage[utf8]{inputenc}
\usepackage[plainpages=false,pdfpagelabels=true]{hyperref}
\usepackage{amssymb,amsthm,wasysym}
\usepackage{slashed}
\usepackage[margin=1in]{geometry}

\numberwithin{equation}{section}

\newtheorem{Satz}{Theorem}[section]
\newtheorem{Prop}[Satz]{Proposition}
\newtheorem{Lem}[Satz]{Lemma}

\newcommand{\cR}{{\mathcal R}}
\theoremstyle{definition}
\newtheorem{Dfn}[Satz]{Definition}
\newtheorem{Bem}[Satz]{Remark}
\newtheorem{Bsp}[Satz]{Example}

\allowdisplaybreaks[1]

\renewcommand{\epsilon}{\varepsilon}

\newcommand{\R}{\ensuremath{\mathbb{R}}}
\newcommand{\C}{\ensuremath{\mathbb{C}}}
\newcommand{\D}{\slashed{D}}
\newcommand{\p}{\slashed{\partial}}
\newcommand{\sff}{\mathrm{I\!I}}

\title{Some aspects of Dirac-harmonic maps with curvature term}
\author{Volker Branding}
\date{\today}
\address{TU Wien\\
Institut für diskrete Mathematik und Geometrie\\
Wiedner Hauptstraße 8–10, A-1040 Wien}
\email[]{volker@geometrie.tuwien.ac.at}
\subjclass[2010]{53C27, 58E20, 35J61}
\keywords{Dirac-harmonic map with curvature term, regularity, vanishing theorem}
\begin{document}

\begin{abstract}
We study several geometric and analytic aspects of Dirac-harmonic maps with curvature term
from closed Riemannian surfaces.
\end{abstract} 

\maketitle

\section{Introduction and Results}
Dirac-harmonic maps arise as critical points of part of the nonlinear \(\sigma\)-model studied in
quantum field theory \cite{MR2262709}. They form a pair of a map from a Riemann surface to a Riemannian manifold and a vector spinor.
The equations for Dirac-harmonic maps couple the harmonic map equation to spinor fields.
As limiting cases both harmonic maps and harmonic spinors can be obtained. 
Moreover, Dirac-harmonic maps belong to the class of conformally invariant variational problems and thus have a lot of nice properties.

Many important results for Dirac-harmonic maps have already been established.
This includes several analytical results \cite{MR2176464}, \cite{MR2544729}, \cite{springerlink:10.1007/s00526-012-0512-5}, \cite{MR2267756}
and an existence result for uncoupled solutions \cite{ammannginoux}.
The boundary value problem for Dirac-harmonic maps is discussed in \cite{bvp}.
A heat-flow approach to Dirac-harmonic maps was studied recently in \cite{regularized}, \cite{regularized-surface}.

However, the full nonlinear supersymmetric \(\sigma\)-model in physics contains additional terms,
which are not captured by the analysis of Dirac-harmonic maps, see \cite{MR1707282}, \cite{Callan:1989nz}
for the physical background.
Taking into account and additional two-form in the action functional the resulting equations were studied in \cite{magnetic},
Dirac-harmonic maps to target spaces with torsion are analyzed in \cite{torsion}.

In this article we focus on Dirac-harmonic maps coupled to a curvature term, which were
introduced in \cite{MR2370260}. This set of equations also has an interesting limit.
In the case of the map part being trivial it reduces to a nonlinear Dirac equation,
which was studied in \cite{MR2390834} and \cite{MR2661574}.
Moreover, it should be noted that this equation also appears in the context 
of the spinorial representation of surfaces in \(\R^3\) \cite{MR1653146}
and the Thirring model in quantum field theory  \cite{Thirring195891}.

In the general case Dirac-harmonic maps with curvature term are more complicated then Dirac-harmonic
maps since they consist of a pair of two non-linear equations.

The aim of this article is to establish some basic results for Dirac-harmonic maps with curvature
term, in particular the regularity of weak solutions.

This article is organized as follows. In Section 2 we recall the notion of Dirac-harmonic 
maps with curvature term. Section 3 discusses geometric and Section 4 
analytical properties of Dirac-harmonic maps with curvature term.

\section{Dirac-harmonic maps with curvature term}
Let us now describe the setup in more detail. For a map \(\phi\colon M\to N\)
we study its differential \(d\phi\in\Gamma(T^*M\otimes\phi^{-1}TN)\), integrating the square of its
norm leads to the usual harmonic energy.
We assume that \((M,h)\) is a closed Riemannian spin surface with spinor bundle \(\Sigma M\), 
for more details about spin geometry see the book \cite{MR1031992}. Moreover, let \((N,g)\) be another closed Riemannian
manifold. Together with the pullback bundle \(\phi^{-1}TN\) 
we consider the twisted bundle \(\Sigma M\otimes\phi^{-1}TN\). 
The induced connection on this bundle will be denoted by \(\tilde{\nabla}\). 
Sections \(\psi\in\Gamma(\Sigma M\otimes\phi^{-1}TN)\) in this bundle  are called \emph{vector spinors} and the natural operator acting on them 
is the twisted Dirac operator, denoted by \(\D\).  
It is an elliptic, first order operator, which is self-adjoint with respect to the \(L^2\)-norm.
More precisely, the twisted Dirac operator is given by \(\D=e_\alpha\cdot\tilde{\nabla}_{e_\alpha}\), where \(\{e_\alpha\}\) is
an orthonormal basis of \(TM\) and \(\cdot\) denotes Clifford multiplication. 
We are using the Einstein summation convention, that is we sum over repeated indices.
Clifford multiplication is skew-symmetric, namely
\[
\langle\chi,X\cdot\xi\rangle_{\Sigma M}=-\langle X\cdot\chi,\xi\rangle_{\Sigma M}
\]
for all \(\chi,\xi\in\Gamma(\Sigma M)\) and all \(X\in TM\).
In addition, the Clifford relations
\[
X\cdot Y+Y\cdot X=-2h(X,Y)
\]
hold for all \(X,Y\in TM\).

We use Greek letters for indices on \(M\) and Latin letters for indices on \(N\).
In terms of local coordinates \(y^i\) the vector spinor \(\psi\) can be written as \(\psi=\psi^i\otimes\frac{\partial}{\partial y^i}\), 
and thus the twisted Dirac operator \(\D\) is locally given by
\[
\D\psi=\big(\p\psi^i+\Gamma^i_{jk}\nabla\phi^j\cdot\psi^k\big)\otimes\frac{\partial}{\partial y^i}.
\]
Here, \(\p\colon\Gamma(\Sigma M)\to\Gamma(\Sigma M)\) denotes the usual Dirac operator.
Throughout this article we study the functional
\begin{equation}
\label{energy-functional}
E_c(\phi,\psi)=\frac{1}{2}\int_M|d\phi|^2+\langle\psi,\D\psi\rangle-\frac{1}{6}\langle R^N(\psi,\psi)\psi,\psi\rangle.
\end{equation}
The factor of \(1/6\) in front of the curvature term is required by supersymmetry.
Here, the indices are contracted as 
\[
\langle R^N(\psi,\psi)\psi,\psi\rangle=R_{ijkl}\langle\psi^i,\psi^k\rangle\langle\psi^j,\psi^l\rangle,
\]
winch ensures that the functional is real valued.
We briefly recall the derivation of the critical points, see \cite{MR2370260}, Section II.
\begin{Prop}
The critical points of the energy functional (\ref{energy-functional}) are given by
\begin{align}
\label{euler-lagrange-phi}\tau(\phi)=&\cR(\phi,\psi)-\tilde{\cR}(\psi), \\
\label{euler-lagrange-psi}\D\psi=&\frac{1}{3}R^N(\psi,\psi)\psi,
\end{align}
where \(\tau(\phi)\) is the tension field of the map \(\phi\) and the curvature terms \(\cR(\phi,\psi)\) 
and \(\tilde{\cR}(\psi)\) are explicitly given by
\begin{align*}
\cR(\phi,\psi)=&\frac{1}{2}R^N(\psi,e_\alpha\cdot\psi)d\phi(e_\alpha),\\
\tilde{\cR}(\psi)=&\frac{1}{12}\langle(\nabla R^N)^\sharp(\psi,\psi)\psi,\psi\rangle.
\end{align*}
Here, \(R^N\) denotes the curvature tensor on \(N\) and \(\sharp\colon\phi^{-1}T^\ast N\to\phi^{-1}TN\) represents the musical isomorphism.
\end{Prop}

\begin{proof}
Consider a variation of \((\phi,\psi)\), that is \(\phi_t\colon(-\epsilon,\epsilon)\times M\to N\) and
\(\psi_t\colon(-\epsilon,\epsilon)\times M\to\Sigma M\otimes\phi_t^{-1}TN\) satisfying 
\((\frac{\partial\phi_t}{\partial t},\frac{\tilde{\nabla}\psi_t}{\partial t})\big|_{t=0}=(\eta,\xi)\).
Here, \(\eta\in\Gamma(\phi^{-1}TN)\) and \(\xi\in\Gamma(\Sigma M\otimes\phi^{-1}TN)\).
It is well known that
\[
\frac{\partial}{\partial t}\big|_{t=0}\frac{1}{2}\int_M|d\phi_t|^2=-\int_M\langle\tau(\phi),\eta\rangle.
\]
The variation of the twisted Dirac-energy was already calculated in \cite{MR2262709}, Prop.2.1, and yields
\[
\frac{\partial}{\partial t}\big|_{t=0}\frac{1}{2}\int_M\langle\psi_t,\D\psi_t\rangle=\int_M\langle\eta,\cR(\phi,\psi)\rangle+\langle\D\psi,\xi\rangle.
\]
Finally, we calculate the variation of the curvature term
\begin{align*}
\frac{\partial}{\partial t}\big|_{t=0}\frac{1}{12}\int_M\langle R^N(\psi_t,\psi_t)\psi_t,\psi_t\rangle =& \frac{1}{12}\int_M\langle (\nabla_{d\phi_t(\partial_t)}R^N)(\psi_t,\psi_t)\psi_t,\psi_t\rangle\big|_{t=0} 
+\frac{4}{12}\int_M\langle R^N(\psi,\psi)\psi,\xi\rangle.
\end{align*}
Concerning the first term on the right hand side, we calculate
\begin{align*}
\langle(\nabla_{d\phi_t(\partial_t)}R^N)(\psi_t,\psi_t)\psi_t,\psi_t\rangle_{\Sigma M\otimes\phi_t^{-1}TN}\big|_{t=0}
&=\big\langle\langle(\nabla R^N)(\psi_t,\psi_t)\psi_t,\psi_t\rangle_{\Sigma M\otimes\phi_t^{-1}TN},\frac{\partial\phi_t}{\partial t}\big\rangle_{\phi_t^{-1}TN}\big|_{t=0} \\
&:=\langle\tilde{R}(\psi),\eta\rangle_{\phi^{-1}TN}
\end{align*}
and adding up the three terms completes the proof.
\end{proof}

In terms of local coordinates the Euler Lagrange equations acquire the form:
\begin{align}
\label{phi-local-coordinates}
\tau^m(\phi)=&\frac{1}{2}R^m_{~lij}\langle\psi^i,\nabla\phi^l\cdot\psi^j\rangle
-\frac{1}{12}g^{mp}(\nabla_pR_{ijkl})\langle\psi^i,\psi^k\rangle\langle\psi^j,\psi^l\rangle, \\
\label{psi-local-coordinates}
(\D\psi)^m=&\frac{1}{3}R^m_{~jkl}\langle\psi^j,\psi^l\rangle\psi^k.
\end{align}
Solutions of the system \eqref{phi-local-coordinates}, \eqref{psi-local-coordinates} are called
\emph{Dirac-harmonic maps with curvature term}.

It would be nice to have some explicit solutions of the Euler-Lagrange equations.
We cannot give an example for a two-dimensional domain, however in the one-dimensional case we have:
\begin{Bsp}
Let \(M=S^1\) and \(\gamma\colon S^1\to N\) be a curve.
Clifford multiplication on \(S^1\) is given by multiplication with the imaginary unit \(i\). 
On \(S^1\) there exist two spin structures, for the trivial spin structure spinors can be identified 
with periodic function, whereas for the non-trivial spin structure spinors are represented by 
anti-periodic functions. Then it can be checked by a direct calculation that
a solution of \eqref{phi-local-coordinates}, \eqref{psi-local-coordinates} is given by
\begin{align}
\tau(\gamma)=0,\qquad \psi=i\chi\otimes\gamma'.
\end{align}
Here, \(\chi\) is a constant spinor and \('\) denotes differentiation with respect to
the curve parameter. This solution is uncoupled in the sense that all terms in the
Euler-Lagrange equations vanish independently.
Moreover, a constant spinor only exists for the trivial spin structure.
\end{Bsp}

\section{Geometric Aspects of Dirac-harmonic Maps with curvature term}
In this section we discuss geometric aspects of Dirac-harmonic maps with curvature term from surfaces.
\begin{Lem}
For a two dimensional domain the functional \(E_c(\phi,\psi)\) is conformally invariant.
\end{Lem}
\begin{proof}
It is well-known that all the three terms in \(E_c(\phi,\psi)\) are invariant
under conformal transformations. For a proof the reader may take a look at Lemma 3.1
in \cite{MR2262709}.
\end{proof}
The energy-momentum tensor for the functional \(E_c(\phi,\psi)\) is given by
\begin{equation}
T_{\alpha\beta}=2\langle d\phi(e_\alpha),d\phi(e_\beta)\rangle-h_{\alpha\beta}|d\phi|^2
+\langle\psi,e_\alpha\cdot\tilde{\nabla}_{e_\beta}\psi\rangle
-\frac{1}{6}h_{\alpha\beta}\langle R^N(\psi,\psi)\psi,\psi\rangle.
\end{equation}
It can be read of from the definition that \(T_{\alpha\beta}\) is traceless, when \((\phi,\psi)\)
is a Dirac-harmonic map with curvature term. In order to see that the energy momentum tensor is 
symmetric, we note that
\[
T_{12}-T_{21}=\langle\psi,e_1\cdot\tilde{\nabla}_{e_2}\psi-e_2\cdot\tilde{\nabla}_{e_1}\psi\rangle
=-\frac{1}{3}\langle\psi,e_2\cdot e_1\cdot R^N(\psi,\psi)\psi\rangle,
\]
where we used that \(\psi\) is a solution of \eqref{euler-lagrange-psi}. Moreover, we have
\begin{align*}
\overline{R_{ijkl}\langle\psi^i,\psi^k\rangle\langle\psi^j,e_1\cdot e_2\cdot\psi^l\rangle}=&R_{ijkl}\langle\psi^k,\psi^i\rangle\langle e_1\cdot e_2\cdot\psi^l,\psi^j\rangle \\
=&-R_{klij}\langle\psi^k,\psi^i\rangle\langle\psi^l,e_1\cdot e_2\cdot\psi^j\rangle \\
=&-R_{ijkl}\langle\psi^i,\psi^k\rangle\langle\psi^j,e_1\cdot e_2\cdot\psi^l\rangle.
\end{align*}
The expression \(T_{12}-T_{21}\) is both purely real and purely imaginary and thus vanishes.
Moreover, a direct calculation gives the following
\begin{Prop}
Let \((\phi,\psi)\) be a Dirac-harmonic map with curvature term. Then the energy momentum tensor
is covariantly conserved, that is
\[
\nabla_{e_\alpha}T_{\alpha\beta}=0 .
\]
\end{Prop}
\begin{proof}
We set
\begin{align*}
C_{\alpha\beta}:=2\langle d\phi(e_\alpha),d\phi(e_\beta)\rangle-h_{\alpha\beta}|d\phi|^2,~~~~
D_{\alpha\beta}:=\langle\psi,e_\alpha\cdot\tilde{\nabla}_{e_\beta}\psi\rangle,~~~~
E_{\alpha\beta}:=-\frac{1}{6}h_{\alpha\beta}\langle R^N(\psi,\psi)\psi,\psi\rangle.
\end{align*}
By a direct calculation, we find
\begin{align*}
\nabla_{e_\alpha}C_{\alpha\beta}=2\langle\tau(\phi),d\phi(e_\beta)\rangle 
=2\langle\cR(\phi,\psi),d\phi(e_\beta)\rangle-2\langle\tilde{\cR}(\psi),d\phi(e_\beta)\rangle.
\end{align*}
Again, calculating directly, we get
\begin{align*}
\nabla_{e_\alpha}D_{\alpha\beta}
=\langle\tilde{\nabla}_{e_\alpha}\psi,e_\alpha\cdot\tilde{\nabla}_{e_\beta}\psi\rangle +\langle\psi,\D\tilde{\nabla}_{e_\beta}\psi\rangle 
=-\langle\D\psi,\tilde{\nabla}_{e_\beta}\psi\rangle +\langle\psi,\D\tilde{\nabla}_{e_\beta}\psi\rangle.
\end{align*}
On the other hand, we find
\begin{align*}
\langle\psi,\D\tilde{\nabla}_{e_\beta}\psi\rangle=&\langle\psi,\tilde{\nabla}_{e_\beta}\D\psi\rangle
+\underbrace{\langle\psi,e_\alpha\cdot R^{\Sigma M}(e_\alpha,e_\beta)\psi\rangle}_{=\frac{1}{2}\langle\psi,\operatorname{Ric}(e_\beta)\cdot\psi\rangle=0} 
+\langle\psi,e_\alpha\cdot R^{N}(d\phi(e_\alpha),d\phi(e_\beta))\psi\rangle \\
=&\langle\psi,\tilde{\nabla}_{e_\beta}\D\psi\rangle-2\langle\cR(\phi,\psi),d\phi(e_\beta)\rangle.
\end{align*}
Finally, we calculate
\begin{align*}
\nabla_{e_\alpha}E_{\alpha\beta}=-\nabla_{e_\alpha}\frac{1}{6}h_{\alpha\beta}\langle R^N(\psi,\psi)\psi,\psi\rangle
=-\frac{1}{2}\nabla_{e_\beta}\langle\psi,\D\psi\rangle
=-\frac{1}{2}\langle\tilde{\nabla}_{e_\beta}\psi,\D\psi\rangle-\frac{1}{2}\langle\psi,\tilde{\nabla}_{e_\beta}\D\psi\rangle.
\end{align*}
Adding up the terms, we find
\begin{align*}
\nabla_{e_\alpha}T_{\alpha\beta}=-2\langle\tilde{\cR}(\psi),d\phi(e_\beta)\rangle
+\frac{1}{2}\langle\psi,\tilde{\nabla}_{e_\beta}\D\psi\rangle-\frac{3}{2}\langle\D\psi,\tilde{\nabla}_{e_\beta}\psi\rangle.
\end{align*}
Using the equation Euler Lagrange equation for \(\psi\) again, we may deduce
\begin{align*}
\langle\psi,\tilde{\nabla}_{e_\beta}\D\psi\rangle&=\frac{1}{3}\langle\psi,\tilde{\nabla}_{e_\beta}(R^N(\psi,\psi)\psi)\rangle\\
&=4\langle\tilde{\cR}(\psi),d\phi(e_\beta)\rangle+\langle R^N(\psi,\psi)\psi,\tilde{\nabla}_{e_\beta}\psi\rangle \\
&=4\langle\tilde{\cR}(\psi),d\phi(e_\beta)\rangle+3\langle \D\psi,\tilde{\nabla}_{e_\beta}\psi\rangle 
\end{align*}
and this proves the assertion.
\end{proof}
On a small domain \(\tilde{M}\) of \(M\) we choose a local isothermal parameter \(z=x+iy\) and set
\begin{align}
\label{hopf-differential}
T(z)dz^2=&(|\phi_x|^2-|\phi_y|^2-2i\langle\phi_x,\phi_y\rangle 
+\langle\psi,\partial_x\cdot\tilde{\nabla}_{\partial_x}\psi\rangle-i\langle\psi,\partial_x\cdot\tilde{\nabla}_{\partial_y}\psi\rangle 
-\frac{1}{3}\langle R^N(\psi,\psi)\psi,\psi\rangle)dz^2
\end{align}
with \(\partial_x=\frac{\partial}{\partial x}\) and \(\partial_y=\frac{\partial}{\partial y}\).
\begin{Prop}
\label{prop-energy-momentum-preserved}
The quadratic differential \(T(z)dz^2\) is holomorphic.
\end{Prop}
\begin{proof}
This follows directly from the last Lemma.
\end{proof}

In addition, we also have 
\begin{Lem}[Bochner formula]
Let \(\psi\) be  a solution of \eqref{euler-lagrange-psi}. Then the following formula holds
\begin{equation}
\Delta\frac{1}{2}|\psi|^2=|\tilde{\nabla}\psi|^2+\frac{R}{2}|\psi|^2+\frac{1}{2}\langle e_\alpha\cdot e_\beta\cdot R^N(d\phi(e_\alpha),d\phi(e_\beta))\psi,\psi\rangle
+\frac{1}{9}|R^N(\psi,\psi)\psi|^2,
\end{equation}
where \(R\) denotes the scalar curvature of \(M\).
\end{Lem}
\begin{proof}
Recall the Weitzenböck formula for the twisted Dirac-operator \(\D\)
\begin{equation}
\label{weitzenboeck-diracoperator}
\D^2\psi=-\tilde{\Delta}\psi+\frac{1}{4}R\psi +\frac{1}{2}e_\alpha\cdot e_\beta\cdot R^N(d\phi(e_\alpha),d\phi(e_\beta))\psi.
\end{equation}
This formula can be deduced from the general Weitzenböck formula for twisted Dirac operators,
see \cite{MR1031992}, p.\ 164, Theorem 8.17.	
Moreover, applying \(\D\) to \eqref{euler-lagrange-psi} we find
\[
\D^2\psi=\frac{1}{3}\big(\nabla\big(R^N(\psi,\psi)\psi\big)\big)\cdot\psi+\frac{1}{3}R^N(\psi,\psi)\D\psi
\]
yielding
\[
-\langle\psi,\D^2\psi\rangle=-\frac{1}{3}\langle\big(\nabla\big(R^N(\psi,\psi)\psi\big)\big)\cdot\psi,\psi\rangle-\frac{1}{9}\langle R^N(\psi,\psi)R^N(\psi,\psi)\psi,\psi\rangle.
\]
The first term on the right hand side vanishes since it is both purely imaginary and purely real.
The claim then follows by a direct calculation.
\end{proof}

As a next step we study the Euler-Lagrange equations for the case that the target manifold \(N\) is isometrically embedded in some \(\R^q\)
by the Nash embedding theorem. Then, we have that \(\phi\colon M\to\R^q\) with \(\phi(x)\in N\). The vector spinor \(\psi\) becomes
a vector of usual spinors \(\psi^1,\psi^2,\ldots,\psi^q\), more precisely \(\psi\in\Gamma(\Sigma M\otimes T\R^q)\).
The condition that \(\psi\) is along the map \(\phi\) is now encoded as
\[
\sum_{i=1}^q\nu^i\psi^i=0\qquad \text{for any normal vector }\nu \text{ at } \phi(x).
\]

\begin{Lem}
Suppose that \(N\subset\R^q\). Then the Euler-Lagrange equations acquire the form
\begin{align}
\label{phirq}-\Delta\phi=&\sff(d\phi,d\phi)+P(\sff(e_\alpha\cdot\psi,d\phi(e_\alpha)),\psi)-G(\psi), \\
\label{psirq}\p\psi=&\sff(d\phi(e_\alpha),e_\alpha\cdot\psi)+F(\psi,\psi)\psi
\end{align}
with the terms
\begin{align*}
G(\psi)&=\frac{1}{6}\langle\nabla_p \sff(\partial_{y^i},\partial_{y^k}),\sff(\partial_{y^j},\partial_{y^l})\rangle-
\langle \nabla_p\sff(\partial_{y^i},\partial_{y^l}),\sff(\partial_{y^j},\partial_{y^k})\rangle)\langle\psi^i,\psi^k\rangle\langle\psi^j,\psi^l\rangle, \\
F(\psi,\psi)\psi&=\frac{1}{3}(P(\sff(\partial_{y^k},\partial_{y^j}),\partial_{y^l})-P(\sff(\partial_{y^l},\partial_{y^j}),\partial_{y^k}))\langle\psi^j,\psi^l\rangle\psi^k.
\end{align*}
Here, \(\sff\) denotes the second fundamental form of \(\phi\) in \(\R^q\) and \(P\) the shape operator.
\end{Lem}
We have written up the terms \(F(\psi,\psi)\psi\) and \(G(\psi)\) in index notation
to highlight the order in which the spinors are contracted.
\begin{proof}
We only discuss the terms arising from the curvature term of \eqref{energy-functional}.
For the other terms see \cite{MR2176464}, p.64-65.
Thus, suppose that \(N\subset\R^q\). Combing both Gauss and Weingarten equation, we get for \(X,Y,Z,W\in\Gamma(TN)\)
\begin{align*}
\langle R^N(X,Y)Z,W\rangle=&\langle\sff(X,Z),\sff(Y,W)\rangle-\langle\sff(Y,Z),\sff(X,W)\rangle \\ 
=&\langle P(\sff(X,Z),Y),W\rangle-\langle P(\sff(Y,Z),X),W\rangle.
\end{align*}
Since \(W\) is arbitrary we may follow
\[
R^N(X,Y)Z=P(\sff(X,Z),Y)-P(\sff(Y,Z),X).
\]
This establishes the equation for \(\psi\). Concerning the equation for \(\phi\),
we make use of the Gauss equation again. Assume that \(X,Y,Z,W\in\Gamma(TN)\) are parallel vector fields.
By differentiation we obtain
\begin{align*}
\langle (\nabla R^N)(X,Y)Z,W\rangle=&\langle(\nabla\sff)(X,Z),\sff(Y,W)\rangle +\langle\sff(X,Z),(\nabla\sff)(Y,W)\rangle \\
&-\langle(\nabla\sff)(Y,Z),\sff(X,W)\rangle -\langle\sff(Y,Z),(\nabla\sff)(X,W)\rangle.
\end{align*}
Furthermore, contracting the indices, we find
\begin{align*}
\nabla_p R_{ijkl}\langle\psi^i,\psi^k\rangle\langle\psi^j,\psi^l\rangle
=&2\langle\nabla_p \sff_{ik},\sff_{jl}\rangle-\langle \nabla_p\sff_{il},\sff_{jk}\rangle)\operatorname{Re}(\langle\psi^i,\psi^k\rangle\langle\psi^j,\psi^l\rangle),
\end{align*}
where \(\sff_{ij}=\sff(\partial_{y^i},\partial_{y^j})\).
\end{proof}

Making use of \eqref{psirq} we may also follow
\begin{Lem}
The zero-set of solutions of \eqref{psirq} is discrete.
\end{Lem}
\begin{proof}
We can think of \eqref{psirq} as an equation of the form
\(
\p\psi=V(\psi).
\)
This is an semielliptic equation of first order. Moreover, \(V\) respects the zero
section, that is \(V(0)=0\). The claim then follows from the main theorem in \cite{MR1714341}.
\end{proof}

\section{Analytic Aspects of Dirac-harmonic maps with curvature term}
This section is devoted to analytic aspects of Dirac-harmonic maps with curvature term.
\subsection{Regularity of solutions}
To study the regularity of Dirac-harmonic maps with curvature term 
from closed Riemannian spin surfaces we define
\begin{align*}
\chi(M,N):=\{(\phi,\psi)\in W^{1,2}(M,N)\times W^{1,\frac{4}{3}}(M,\Sigma M\otimes\phi^{-1}TN) 
\text{ with } (\ref{phirq}) \text{ and } (\ref{psirq}) \text{ a.e.}\}.
\end{align*}
\begin{Dfn}
A pair \((\phi,\psi)\in\chi(M,N)\) is called \emph{weak Dirac-harmonic map with curvature term} from \(M\) to \(N\) if and only if the pair \((\phi,\psi)\) solves
\eqref{phirq}, \eqref{psirq} in a distributional sense.
\end{Dfn}
First, we improve the regularity of \(\psi\) with the help of the tools presented in \cite{MR2661574}.
To this end we recall the definition of Morrey spaces.
\begin{Dfn}
Assume that \(1\leq p\leq n, 0<\lambda\leq n\) and let \(U\subset \R^n\) be a domain.
Then the Morrey space \(M^{p,\lambda}(U)\) is defined by
\[
M^{p,\lambda}(U):=\big\{f\in L^p_{loc}(U)\mid |f|_{M^{p,\lambda}(U)}<\infty\big\},
\]
where 
\[
|f|^p_{M^{p,\lambda}(U)}:=\sup_{r>0}\big\{r^{\lambda-n}\int_{D_r}|f|^p\mid D_r\subset U\big\}.
\]
\end{Dfn}
Note that for \(1\leq p\leq n\) we have \(M^{p,\lambda}(U)\subset L^p(U)\)
and \(M^{p,n}(U)=L^p(U)\).
Moreover, we also recall the Riesz-potential of order \(1\)
\[
I_1(f)(x)=\int_{\R^2}\frac{|f(y)|}{|x-y|},\qquad f\colon\R^2\to\R.
\]
In addition, we want to apply Adam's inequality on Morrey spaces (see \cite{MR0458158}, Theorem 3.1)
\begin{equation}
\label{adams-riesz}
|I_1(f)|_{M^{\frac{\lambda q}{\lambda-q},\lambda}(\R^n)}\leq C|f|_{M^{q,\lambda}(\R^n)},\qquad 1\leq q\leq\lambda<n.
\end{equation}
We now derive an estimate for the spinor \(\psi\), therefore we set
\[
A(\phi,\psi)\psi=\sff(e_\alpha\cdot\psi,d\phi(e_\alpha))+F(\psi,\psi)\psi.
\]
In the following \(D_r\subset \R^2\) will denote the disc with radius \(r\) with \(r\leq 1\).
The next Lemma follows Lemma 2.2 in \cite{MR2661574}, we thus do not give all the details,
see also \cite{1306.4260}, Theorem 4.1.
Note that we may trivialize the spinor bundle \(\Sigma M\) on the disc \(D_r\) by two complex functions.
\begin{Prop}
\label{lemma-spinor-morey}
Let \(D_1\subset\R^2\) and \(\psi\in L^4(D_1,\C^2\otimes\R^q), A(\phi,\psi)\in L^2(D_1,gl(q)\otimes\R^2)\).
Suppose that \(\psi\) weakly solves
\[
\p\psi=A(\phi,\psi)\psi
\]
on the disc \(D_1\). If
\[
|A(\phi,\psi)|_{L^2(D_1)}\leq\epsilon
\]
then for any \(4<p<\infty\) the following estimate holds
\begin{equation}
\label{nabla-psi-lq}
|\psi|_{L^p(D_\frac{1}{16})}\leq C|\psi|_{L^4(D_1)}.
\end{equation}
\end{Prop}

\begin{proof}
Applying the usual Dirac operator to the first order equation \eqref{psirq} yields
\[
-\Delta\psi=\p^2\psi=\p(A(\phi,\psi)\psi)
\]
in a distributional sense, where we applied the Schrödinger-Lichnerowicz formula and
used the fact that the disc is flat.
We set \(V(\phi,\psi):=N[\p(A(\phi,\psi)\psi)]\), where \(N\) denotes the Newtonian potential and get
\(|V(\phi,\psi)|\leq I_1[|A(\phi,\psi)\psi|]\).
Moreover, we find
\begin{align}
|V(\phi,\psi)|_{L^4(\R^2)}&\leq I_1[|A(\phi,\psi)\psi|]_{L^4(\R^2)} \\
\nonumber &\leq C|A(\phi,\psi)\psi|_{L^\frac{4}{3}(\R^2)}\\
\nonumber &\leq C|A(\phi,\psi)|_{L^2(\R^2)}|\psi|_{L^4(\R^2)},
\end{align}
where we applied \eqref{adams-riesz}. This estimate can be localized to the disc with the help of a cutoff function.
As a next step, we set \(Z(\phi,\psi):=V(\phi,\psi)-\psi\) and thus \(\Delta Z(\phi,\psi)=0\).
Using an estimate for harmonic functions (see \cite{MR1913803}, Lemma 3.3.12) we get
\begin{align}
|Z(\phi,\psi)|_{L^4(D_r)}&\leq r^\frac{1}{2}|Z(\phi,\psi)|_{L^4(D_1)}\\
\nonumber &\leq r^\frac{1}{2}(C|A(\phi,\psi)|_{L^2(D_1)}|\psi|_{L^4(D_1)}+|\psi|_{L^4(D_1)}).
\end{align}
Moreover, we may estimate
\begin{align}
|\psi|_{L^4(D_r)}&\leq |V(\phi,\psi)|_{L^4(D_r)}+|Z(\phi,\psi)|_{L^4(D_r)} \\
\nonumber&\leq C\epsilon^2|\psi|_{L^4(D_1)}+Cr^\frac{1}{2}|Z(\phi,\psi)|_{L^4(D_1)} \\
\nonumber&\leq C\epsilon^2|\psi|_{L^4(D_1)}+Cr^\frac{1}{2}(|\psi|_{L^4(D_1)}+|V(\phi,\psi)|_{L^4(D_1)}) \\
\nonumber&\leq C(\epsilon^2+r^\frac{1}{2})|\psi|_{L^4(D_1)}.
\end{align}
Performing a suitable rescaling we thus find
\begin{align*}
|\psi|_{L^4(D_\frac{r}{2})}\leq\theta|\psi|_{L^4(D_r)}
\end{align*}
for \(\theta\) chosen appropriately.
By iteration this yields
\[
|\psi|_{L^4(D_r)}\leq r^\frac{\alpha}{2}|\psi|_{L^4(D_1)}\qquad\text{ for}\qquad\text0<\alpha<\frac{1}{3}, 0<r<\frac{1}{4}
\]
and we may follow that
\[
|\psi|_{M^{4,2-\nu}(D_\frac{1}{8})}\leq C|\psi|_{L^4(D_1)} 
\]
for \(\nu\in(0,1)\).
At this point we apply Adam's inequality \eqref{adams-riesz} again with \(q=4/3\) and \(\lambda=2-\nu\).
Repeating the procedure from above we find
\[
|\psi|_{M^{\frac{4(2-\nu)}{2-3\nu},2-\nu}(D_\frac{1}{16})}\leq C|\psi|_{L^4(D_1)}.
\]
Now, we let \(\nu\to 2/3\) and may follow
\begin{equation}
|\psi|_{L^p(D_\frac{1}{16})}\leq C|\psi|_{L^4(D_1)} 
\end{equation}
for any \(p>4\), which completes the proof.
For more details see the end of the proof of Lemma 2.2 in \cite{MR2661574}.
\end{proof}
Note that \(|A(\phi,\psi)|_{L^2}\leq C\) for \((\phi,\psi)\in\chi(M,N)\).

As a second step, we turn to the Euler-Lagrange equation for \(\phi\).
Following \cite{bvp}, p.7, we can rewrite the equation \eqref{phirq} as follows
\begin{align}
\label{phi-antisym}
-\Delta\phi^m=&\phi^i_\alpha\phi^j_\alpha\bigg(\frac{\partial\nu_l^i}{\partial y^j}\nu_l^m-\frac{\partial\nu_l^m}{\partial y^j}\nu_l^i\bigg)-G^m(\psi) \\
\nonumber&+\phi^i_\alpha\langle\psi^k,e_\alpha\cdot\psi^j\rangle 
\bigg(\bigg(\frac{\partial\nu_l}{\partial y^j}\bigg)^{\top,i}\bigg(\frac{\partial\nu_l}{\partial y^k}\bigg)^{\top,m}-\bigg(\frac{\partial\nu_l}{\partial y^k}\bigg)^{\top,i}\bigg(\frac{\partial\nu_l}{\partial y^j}\bigg)^{\top,m}\bigg).
\end{align}
Here, \(\top\) denotes the projection map \(\top\colon\R^q\to T_yN\).
Moreover, \(\nu_l,l=n+1,\ldots,q\) is an orthonormal frame field for the normal bundle \(T^\perp N\)
and \(\phi_\alpha\) represents the derivative of \(\phi\) with respect to \(e_\alpha\).

Thus, the equation for \(\phi\) can be written in the form
\begin{equation}
-\Delta\phi^m=A^m_{~~i}\cdot\nabla\phi^i+G^m(\psi)\qquad \text{        with       }    A^m_{~~i}=-A^i_{~~m}.
\end{equation}
Here,
\[
A^m_{~i}=\begin{pmatrix}
f^m_{~i} \\
g^m_{~i}
\end{pmatrix},
\qquad
i,m=1,2,\ldots,q
\]
with
\begin{align*}
f^m_{~i}:=&\left(\frac{\partial\nu_l^i}{\partial y^j}\nu_l^m-\frac{\partial\nu_l^m}{\partial y^j}\nu_l^i\right)\phi^j_x+\langle\psi^k,\partial_x\cdot\psi^j\rangle
\bigg(\bigg(\frac{\partial\nu_l}{\partial y^j}\bigg)^{\top,i}\bigg(\frac{\partial\nu_l}{\partial y^k}\bigg)^{\top,m}-\bigg(\frac{\partial\nu_l}{\partial y^k}\bigg)^{\top,i}\bigg(\frac{\partial\nu_l}{\partial y^j}\bigg)^{\top,m}\bigg)
\end{align*}
and we get a similar equation for \(g^m_{~i}\).
The antisymmetry of \(A^m_{~~i}\) will be the important point in the following
regularity analysis.

\begin{Bem}
In the case of a spherical target, one can derive the continuity of the map \(\phi\) by ``classical tools''.
Using the symmetries of the target \(S^n\) it is possible to write the equation for \(\phi\) in such a form
that one may apply the classical Wente Lemma \cite{MR0243467} yielding the continuity of the map \(\phi\)
similar to \cite{MR2176464}, Prop. 2.1.
\end{Bem}

To improve the regularity of the map \(\phi\) we will use the following Theorem from \cite{MR3020100}.
\begin{Satz}
\label{topping-sharp}
Suppose that \(\phi\in W^{1,2}(D_1,\R^q)\) is a weak solution of
\begin{equation}
-\Delta\phi=\Omega\cdot\nabla\phi+f,\qquad f\in L^p(D_1,\R^q),
\end{equation}
where \(\Omega\in L^2(D_1,so(q)\otimes\R^2)\) and \(p\in(1,2)\). Then \(\phi\in W_{loc}^{2,p}(D_1)\).
In particular, if \(f=0\), then \(\phi\in W_{loc}^{2,p}\) for all \(p\in [1,2)\) and \(\phi\in W_{loc}^{1,q}\) for all \(q\in [1,\infty)\).
Moreover, for \(U\subset D_1\), there exist \(\eta_0=\eta_0(p,q)>0\) and \(C=C(p,m,U)<\infty\)
such that if \(|\Omega|_{L^2(D_1)}\leq\eta_0\), then the following estimate holds
\begin{equation}
\label{phi-w2p}
|\phi|_{W^{2,p}(U)}\leq C(|f|_{L^p(D_1)}+|\phi|_{L^1(D_1)}).
\end{equation}
\end{Satz}
Thus, applying Theorem \ref{topping-sharp} to \eqref{phi-antisym} we immediately get that \(\phi\in W_{loc}^{2,p}\) for \(p\in (1,2)\).

\begin{Satz}
Suppose that \((\phi,\psi)\) is a weak Dirac-harmonic map with curvature term from a closed
Riemannian spin surface to a compact Riemannian manifold. Then the pair \((\phi,\psi)\) is smooth.
\end{Satz}
\begin{proof}
This follows from a bootstrap argument using \eqref{phi-w2p} and \eqref{nabla-psi-lq}.
Applying Theorem \ref{topping-sharp} to \eqref{phi-antisym} with \(p=3/2\) we get
that \(\phi\in W^{2,\frac{3}{2}}\) and by the Sobolev embedding theorem this yields \(d\phi\in L^6\).
Now, using \eqref{phirq} and the estimate for the spinor \eqref{nabla-psi-lq}
we have \(\phi\in W^{2,2}\) and thus \(\phi\in W^{1,p}\) for all \(p>1\).
Using elliptic estimates for \eqref{psirq} we may deduce that \(\psi\in W^{1,p}\)
for any \(p>4\) and by Sobolev embedding this yields 
\(\psi\in C^{0,\alpha}\) for some \(\alpha\in (0,1)\).
By iteration we obtain the smoothness of \((\phi,\psi)\).
\end{proof}

\begin{Bem}
In \cite{MR3018163} the authors consider Dirac-harmonic maps coupled to a Ricci-type potential, more precisely,
they study the functional 
\begin{equation}
E(\phi,\psi)=\frac{1}{2}\int_M|d\phi|^2+\langle\psi,\D\psi\rangle+\operatorname{Ric}(\phi)(\psi,\psi)
\end{equation}
and perform a regularity analysis of the critical points using the method of moving frames.
It turns out, that a weak critical point of this functional with \((\phi,\psi)\in\chi(M,N)\) is smooth. 
However, the authors also note that their regularity result does not carry 
over to the potential studied in the case of Dirac-harmonic maps with curvature term. 
\end{Bem}

We close this section by proving two different ``vanishing Theorems''.

\begin{Lem}
Assume that \(N\) has constant curvature and that pair \((\phi,\psi)\) is a smooth solution of
\eqref{phirq} and \eqref{psirq} satisfying
\begin{equation}
\int_M(|d\phi|^2+|\psi|^4)<\epsilon_0
\end{equation}
with \(\epsilon_0\) small enough. Then \(\phi\) is constant and \(\psi\) solves \(\p\psi=F(\psi,\psi)\psi\).
\end{Lem}
This result is similar to Prop.4.3 from \cite{MR2262709}.
\begin{proof}
Since \(N\) has constant curvature the map \(\phi\) now solves
\[
-\Delta\phi=\sff(d\phi,d\phi)+P(\sff(e_\alpha\cdot\psi,d\phi(e_\alpha)),\psi).
\]
Thus, we estimate
\begin{align*}
|\Delta\phi|_{L^\frac{4}{3}}\leq &C\big(\big||d\phi|^2\big|_{L^\frac{4}{3}}+\big||d\phi||\psi|^2\big|_{L^\frac{4}{3}}\big)\\
\leq &C\big(|d\phi|_{L^2}+|\psi|^2_{L^4}\big)|d\phi|_{W^{1,\frac{4}{3}}} \\
\leq&\epsilon_0 C|d\phi|_{W^{1,\frac{4}{3}}}.
\end{align*}
Hence for \(\epsilon_0\) small enough we may conclude that \(\phi\) is constant.
\end{proof}

\begin{Lem}
Assume that the pair \((\phi,\psi)\) is a smooth solution of
\eqref{phirq} and \eqref{psirq} satisfying
\begin{equation}
\int_M(|d\phi|^2+|\psi|^4)<\epsilon_0
\end{equation}
with \(\epsilon_0\) small enough. Moreover, assume that there do not exist harmonic spinors on \(M\).
Then both \(\phi\) and \(\psi\) are trivial.
\end{Lem}

\begin{proof}
As a first step we prove that \(|\psi|_{L^\frac{4}{3}}\leq C|\p\psi|_{L^\frac{4}{3}}\).
By assumption \(0\) is not in the spectrum of \(\p\), thus we may estimate
\[
|\psi|\leq \frac{1}{|\lambda_1|}|\p\psi|,
\]
where \(\lambda_1\) denotes the smallest eigenvalue of the Dirac operator.
Taking the \(L^\frac{4}{3}\) norm then gives the result.
This inequality can also be proven differently, see Lemma 4.1 in \cite{MR2390834}.

Applying elliptic estimates for first order equations (\cite{MR1031992}, Thm. 5.2) to \eqref{psirq} we may estimate
\begin{align*}
|\psi|_{L^4}\leq C|\psi|_{W^{1,\frac{4}{3}}}
\leq &C\big(|\p\psi|_{L^\frac{4}{3}}+|\psi|_{L^\frac{4}{3}}\big) \\
\leq &C\big(\big||\psi|^3\big|_{L^\frac{4}{3}}+\big||\psi||d\phi|\big|_{L^\frac{4}{3}}+|\p\psi|_{L^\frac{4}{3}}\big)\\
\leq &C\big(|\psi|_{L^4}^3+|\psi|_{L^4}|d\phi|_{L^2}\big)\\
\leq &C|\psi|_{L^4}\big(|\psi|^2_{L^4}+|d\phi|_{L^2}\big) \\
\leq &\epsilon_0C|\psi|_{L^4}.
\end{align*}
Hence, for \(\epsilon_0\) small enough \(\psi\) has to vanish. Using \(\psi=0\) and \eqref{phirq} we can estimate
\begin{align*}
|\Delta\phi|_{L^\frac{4}{3}}\leq &C\big||d\phi|^2\big|_{L^\frac{4}{3}}
\leq C|d\phi|^2_{L^2}|d\phi|_{L^4} 
\leq C\epsilon_0|d\phi|_{W^{1,\frac{4}{3}}}
\end{align*}
and thus \(\phi\) has to be constant, whenever \(\epsilon_0\) is small enough.
\end{proof}
The condition that there are no harmonic spinors on \(M\) is satisfied if \(M\)
has positive Euler characteristic. In general, the existence of harmonic spinors
will depend on both metric and chosen spin structure.

\bibliographystyle{plain}
\bibliography{mybib}
\end{document}